\documentclass{amsart}
 \usepackage{latexsym}
\usepackage{amsmath}
\usepackage{amssymb}

 \newtheorem{thm}{Theorem}[section]
 \newtheorem{cor}[thm]{Corollary}

 \theoremstyle{definition}
 
 \theoremstyle{remark}
 
 \newtheorem{problem}{Problem}

\newtheorem{conjt}{Conjecture}

 \numberwithin{equation}{section}

\newcommand{\za}{\zeta}

\newcommand{\ph}{\varphi}
\newcommand{\cph}{C_\varphi}
\newcommand{\hol}{\mathcal{H}ol}
\newcommand{\Dbb}{\mathbb D}
\newcommand{\Tbb}{\mathbb T}

\newcommand{\Cbb}{\mathbb C}

\newcommand{\scw}{\mathcal{S}}
\newcommand{\meal}{m}
\newcommand{\dist}{\rho}

\numberwithin{equation}{section}

\begin{document}

%
%
%
%
%
%
%
%
%

\title[Compact linear combinations of composition operators]
 {Compact linear combinations of composition operators on Hardy spaces}

\author[E. Doubtsov]{Evgueni Doubtsov}

\address{%
Department of Mathematics and Computer Science\\
St.~Petersburg State University\\
Line 14th (Vasilyevsky Island), 29\\
St.~Petersburg 199178\\
Russia}
\email{dubtsov@pdmi.ras.ru}

\author{Dmitry V. Rutsky}
\address{
St.~Petersburg Department of Steklov Mathematical Institute\\
Fontanka 27\\
St.~Petersburg 191023\\
Russia}
\email{rutsky@pdmi.ras.ru}

\thanks{The research in Sections~1 and 3 was supported by Russian Science Foundation (grant No.~24-11-00087),
https://rscf.ru/project/24-11-00087/;
the research in Sections~2 and 4 was supported by Russian Science Foundation (grant No.~23-11-00171),
https://rscf.ru/project/23-11-00171/}




\begin{abstract}
Let $\varphi_j$, $j=1,2, \dots, N$, be holomorphic self-maps of the unit disk $\mathbb{D}$ of $\mathbb{C}$.
We prove that the compactness of a linear combination of the composition operators 
$C_{\varphi_j}: f\mapsto f\circ\varphi_j$
on the Hardy space $H^p(\mathbb{D})$ does not depend on $p$ for $0<p<\infty$.
This answers a conjecture of Choe et al. about the compact differences $C_{\varphi_1} - C_{\varphi_2}$
on $H^p(\mathbb{D})$, $0<p<\infty$.
\end{abstract}

\maketitle
\section{Introduction}\label{s_int}

Let $\hol (\Dbb)$ denote the space of holomorphic functions in the unit disk $\Dbb$ of the
complex plane $\Cbb$, and
let $\scw(\Dbb)$ be the class of all holomorphic self-maps of $\Dbb$.
Given a $\ph \in\scw(\Dbb)$, the composition operator $\cph$ is defined by
\[
\cph f = f \circ \ph
\]
for $f\in \hol(\Dbb)$. Clearly, $\cph$ maps the space $\hol(\Dbb)$ into itself. 

For $0 < p <\infty$, the Hardy space $H^p(\Dbb)$ consists of $f\in \hol(\Dbb)$ such that
\[
\|f\|_{H^p}^p = \sup_{0<r<1} \int_{\Tbb} |f(r\za)|^p d\meal(\za) < \infty,
\]
where $\meal$ is the normalized Lebesgue measure on the unit circle $\Tbb = \partial\Dbb$.

Given a $\ph\in\scw(\Dbb)$, a typical problem is to relate certain operator-theoretic properties of $C_\ph$
to function-theoretic properties of $\ph$. 
Properties of $C_\ph$ on the Hardy spaces are of special interest.
By the classical Littlewood subordination principle,
$C_\ph$ is bounded on $H^p(\Dbb)$, $p>0$, for any symbol $\ph\in\scw(\Dbb)$.
See standard references \cite{CmC95} and \cite{Shj93} for various aspects on the theory of
composition operators on $H^p(\Dbb)$ and on other classical function spaces
in one and several complex variables.

More generally, one considers similar questions for linear combinations of
 $C_{\ph_j}$, $\ph_j\in \scw(\Dbb)$, $j=1, 2,  \dots, N$.
The present paper is motivated by the following problem concerning differences of two composition operators:

\begin{problem}\label{p_prb}
Given a $p>0$, characterize $\ph, \psi\in \scw(\Dbb)$ such that $C_\ph - C_\psi$ is a compact operator on $H^p(\Dbb)$.
\end{problem}

The above problem goes back to J.~H.~Shapiro and C.~Sundberg \cite{SjS90}.
For $p=2$, Problem~\ref{p_prb} was first treated in \cite{CCKY20JFA}.
A more explicit solution 
was then given in \cite{CCKP22TAMS}.
Moreover, it was shown in \cite{NS04} that the solution of Problem~\ref{p_prb}
does not depend on $p$ for $1\le p<\infty$. 
Finally, B.~R.~Choe et al.\ \cite{CCKP22TAMS} formulated the following conjecture.

\begin{conjt}[{\cite[Conjecture~4.3]{CCKP22TAMS}}]\label{cj_conj}
The compactness of $C_\ph - C_\psi$ on $H^p(\Dbb)$ is independent of the parameter $p\in (0,\infty)$.
\end{conjt}

In the present paper, we obtain the following result.

\begin{thm}\label{t_comp_KTp0}
Let $\varphi_j\in \scw(\Dbb)$, $j=1,2, \dots, N$,
and let $T$ denote a linear combination of the composition operators $C_{\ph_j}$.
The compactness of $T: H^p(\Dbb) \to H^p(\Dbb)$ does not depend on $p$ for $0 < p < \infty$.
\end{thm}

In particular, this answers positively Conjecture~\ref{cj_conj}:

\begin{cor}\label{c_conj}
Let $\ph, \psi \in \scw(\Dbb)$.
The compactness of $C_\ph - C_\psi: H^p(\Dbb) \to H^p(\Dbb)$ does not depend on $p$ for $0<p<\infty$.
\end{cor}

\subsection*{Organization of the paper}
In Section~\ref{s_prf_p1}
we formulate and prove Theorem~\ref{t_cmp_2},
a quasi-Banach version of a classical compactness result due to A.~Persson~\cite {Pe64}.
Theorem~\ref{t_cmp_2} implies Corollary~\ref{t_cmpIntrp_Cw}, a one-sided compactness result
for the real interpolation of quasi-Banach Hardy spaces.
We apply Corollary~\ref{t_cmpIntrp_Cw} to prove Theorem~\ref{t_comp_KTp0} in Section~\ref{s_theProof}.
The final Section~\ref{s_final} contains some corollaries from Theorem~\ref{t_comp_KTp0}.

\section{Auxiliary results: compactness and interpolation}\label{s_prf_p1}

\subsection{Real interpolation}
Below we use  
the real interpolation method for quasi-Banach spaces, 
so we first introduce 
related basic notions.

Let $(A_0, A_1)$ be a compatible couple of quasi-Banach spaces.
The interpolation space $(A_0, A_1)_{\theta, q}$,
$0< \theta < 1$, $0 < q \le \infty$, is defined by
\[
(A_0, A_1)_{\theta, q} = \{f\in A_0 + A_1: \|f\|_{\theta, q}<\infty\},
\]
where
\[
\|f\|_{\theta, q}^q =\int_0^\infty \left(t^{-\theta} K(t,f) \right)^q \, \frac{dt}{t}
\]
and
\begin{equation}\label{e_K_df}
\begin{split}
  K(t, f) 
&= K(t, f; A_0, A_1) \\
&= \inf\{\|f_0\|_{A_0} + t\|f_1\|_{A_1} : f= f_0 + f_1,\ f_0\in A_0,\ f_1\in A_1\}
\end{split}
\end{equation}
is the K-functional.
See, for example, \cite{BL76} for further details.

\subsection{A compactness theorem}
The following is a 
generalization of A.~Persson's classical result~\cite {Pe64}.
For the sake of clarity we write it out in detail.
\begin{thm}\label{t_cmp_2}
Let $(E_0, E_1)$ and~$(F_0, F_1)$ be two compatible couples of quasi-Banach spaces, and suppose that
$E$ and~$F$ are quasi-Banach interpolation spaces for these couples, meaning that for any linear operator~$U$ acting
boundedly as~$U : E_0 \to F_0$, $U : E_1 \to F_1$
it follows that it is also bounded as~$U : E \to F$.
Suppose also that the second couple satisfies
the following approximation condition: there exists a sequence of linear operators~$P_j : F_0, F_1 \to F_0 \cap F_1,$ such that 
\begin {equation}
\label {e_H_2}
\tag {H} P_j f \to f \text { in } F_0 \text { for all } f \in F_0 \text { and } \sup_j \|P_j\|_{F_1 \to F_1} < \infty, 
\end{equation}
and the intermediate spaces satisfy the following conditions:
\begin {equation}
\label {e_C_0}
\tag {$C_0$}
\lim_{\varepsilon \to 0} \sup_{{}}
\left\{ \|U\|_{E \to F} : \|U\|_{E_0 \to F_0} \leqslant \varepsilon, \|U\|_{E_1 \to F_1} \leqslant 1 \right\} = 0,
\end{equation}
\begin {equation}
\label {e_C_1}
\tag {$C_1$}
\lim_{t \to \infty} \quad t^{-1} \sup_{f \in B_E} K (t, f; E_0, E_1) = 0,
\end{equation}
where~$B_E$ denotes the (closed) unit ball of~$E$, and
$K (t, f; E_0, E_1)$
is the K-functional defined by \eqref{e_K_df}.
Then any linear operator~$T$ that is compact as~$T : E_0 \to F_0$ and bounded as~$T : E_1 \to F_1$
is also compact as~$T : E \to F$.
\end{thm}
\begin{proof}
We first establish the particular case~$F_0 = F_1 = F$ due to J.~L.~Lions and J.~Peetre (see \cite [Theorem~3.8.1] {BL76}).
Namely, suppose that $T$ is compact as $T: E_0 \to F$
and bounded as $T: E_1 \to F$.
Let~$a_j \in B_E$ and~$\varepsilon > 0$.
By the assumption~\eqref {e_C_1}, for large enough~$t$ there exist decompositions~$a_j = a_{0, j} + a_{1, j}$
satisfying~$\|a_{0, j}\|_{E_0} \leqslant t \varepsilon$ and~$\|a_{1, j}\|_{E_1} \leqslant \varepsilon$. 
By passing to a subsequence, we may assume that
\[
\|T a_{0, j} - T a_{0, j'}\|_F \leqslant \varepsilon
\]
for large enough~$j$ and~$j'$.
On the other hand, 
\[
\|T a_{1, j} - T a_{1, j'}\|_F \leqslant 2 \varepsilon \|T\|_{E_1 \to F}.
\]
Adding these two estimates together yields
\[
\|T a_j - T a_{j'}\|_F \leqslant c (1 + 2 \|T\|_{E_1 \to F}) \varepsilon
\]
 for
large enough~$j$ and~$j'$, where $c$ is a constant in the triangle inequality for~$F$.
A convergent subsequence of~$T a_j$ is then obtained by an application of the diagonal process.

Now, we verify the general case of Theorem~\ref {t_cmp_2}.
Since~$T (B_{E_0})$ is relatively compact in~$F_0$, 
\[
\lim_{j \to \infty} \|P_j T - T\|_{E_0 \to F_0} = 0
\]
 by the assumption~\eqref {e_H_2}.
By the assumption~\eqref {e_C_0},
\[
\lim_{j \to \infty} \|P_j T - T\|_{E \to F} = 0,
\]
that is, $T$ is approximated uniformly by~$P_j T$, and it suffices to verify that~$P_j T : E \to F$ is compact.
Finally,
\[
P_j T : E_0, E_1 \to F_0 \cap F_1 \subset F,
\]
 and~$P_j T : E_0 \to F_0 \cap F_1$ is compact
as a composition of a bounded and a compact operator. 
Thus its compactness is reduced to the particular case~$F_0 = F_1 = F$.
\end{proof}


\begin{cor}\label{t_cmpIntrp_Cw}
Assume that $T: H^{p_j}(\Dbb) \to H^{p_j}(\Dbb)$, $0< p_j < \infty$, $j=0, 1$, is a bounded linear operator such that 
$T:H^{p_1}(\Dbb) \to H^{p_1}(\Dbb)$ is compact. 
Then 
\[
T: (H^{p_0}(\Dbb), H^{p_1}(\Dbb))_{\theta, q}  \to (H^{p_0}(\Dbb), H^{p_1}(\Dbb))_{\theta, q}
\] 
is a compact operator
for all admissible $\theta$ and $q$.
\end{cor}
\begin{proof}
We claim that Theorem~\ref{t_cmp_2} applies.
Indeed, assumptions~\eqref {e_C_0} and~\eqref {e_C_1} are satisfied by all interpolation functors 
of type~$\mathcal C_\theta$, $0 < \theta < 1$, in particular, for $E= (E_0, E_1)_{\theta, q}$
and $F= (F_0, F_1)_{\theta, q}$
(see~\cite[Section~3.5]{BL76}).
Finally, taking~$P_j f (z) = f (r_j z)$ with some sequence~$0 < r_j < 1$, $r_j \to 1$,
we conclude that
assumption~\eqref {e_H_2} is satisfied for the entire range of $H^p(\Dbb)$, $0< p < \infty$.
\end{proof}

\subsection{Remarks on further applications}
Assumptions~\eqref {e_C_0} and~\eqref {e_C_1} are 
also satisfied in many other applications of interest, such as couples of the quasinormed
Orlicz spaces $(L^{\Phi_0}, L^{\Phi_1})$ and the corresponding intermediate Orlicz spaces
\[
L^{\Phi_0} \cap L^{\Phi_1} \subset L^\Phi \subset L^{\Phi_0} + L^{\Phi_1}
\]
such that~$L^{\Phi} \neq L^{\Phi_1}$.
See, for example, \cite {Mit65}.  However, the details in full generality may be somewhat complicated,
so we leave it as a note in passing.

Assumption~\eqref {e_H_2} is easily seen to also be satisfied for the entire range of weighted Bergman spaces 
(excluding~$p = \infty$)
if, as in Corollary~\ref{t_cmpIntrp_Cw}, we take~$P_j f (z) = f (r_j z)$ with some sequence~$0 < r_j < 1$, $r_j \to 1$.
Moreover, the same is true for the Hardy-type spaces for rearrangement invariant quasi-Banach lattices
(see~\cite [Theorem~2.1] {Aleksandrov1981}), which are also stable under the real interpolation
by the results of~\cite {Ki99} (in addition to the stability in the classical case discussed in Section~\ref {s_theProof} below); 
for the $\mathrm {BMO}$-regularity property of these lattices see, for example,
\cite[Proposition~2]{Rut11}.
In particular, Theorem~\ref {t_cmp_2} is also applicable to the Hardy-Orlicz spaces~$H^\Phi(\Dbb)$.

\subsection{Remarks on the compactness theorem}

Theorem~\ref {t_cmp_2} seems to be one of the earliest
(taken from \cite{Pe64})
 and indeed the easiest compactness results for couples of spaces
that is perfectly applicable in many settings.
It should be noted that there exists a large body of work on compactness and interpolation which is also relevant to the research;
for instance, see~\cite{CF89, CEP90, Cw92}.
In particular, the proof of the Hayakawa-Cwikel theorem~\cite[Theorem~1.1]{Cw92},
establishing this result for the standard real interpolation, 
seems to work also in the quasi-Banach setting.
Thus, for these interpolation spaces the assumption~\eqref {e_H_2} is clearly superfluous.
This (\cite[Theorem~1.1]{Cw92}) is indeed the result that one would use for general weighted Hardy spaces,
since for them approximation might not work so easily.
However, whether a compactness theorem holds true for the complex interpolation
in the general case is still a notable open problem in the interpolation theory even for Banach spaces.

\section{Proof of Theorem~\ref{t_comp_KTp0}}\label{s_theProof}
Suppose that $T : H^{p}(\Dbb)\to H^{p}(\Dbb)$ is 
compact with some~$p = p_0\in (0, \infty)$.
We
need to verify
that 
$T: H^{p}(\Dbb)\to H^{p}(\Dbb)$ is also compact for all $0 < p < \infty$.

Fix $p$ and $p_1$ such that $p_1 > p > p_0$ or $p_1 < p < p_0$.
Define $\theta\in (0,1)$ by the following identity: 
\begin{equation}\label{e_theta}
\frac{1}{p} = \frac{1-\theta}{p_0} + \frac{\theta}{p_1}.
\end{equation}
It is well known that \eqref{e_theta} implies
$(H^{p_0}(\Dbb), H^{p_1}(\Dbb))_{\theta, p} = H^{p}(\Dbb)$.
This result essentially goes back to R.~Salem and A.~Zygmund \cite{SZ48} 
(see \cite[Chapter~1, Exercise~14]{BL76}); for a modern treatment and appropriate references see \cite{Ki99}.
Next, recall that $T: H^{p_1}(\Dbb)\to H^{p_1}(\Dbb)$ is a bounded operator
by the Littlewood subordination principle. 
Thus, applying Corollary~\ref{t_cmpIntrp_Cw} with $p_0$, $p_1$ and $q=p$, we conclude that  
$T: H^{p}(\Dbb)\to H^{p}(\Dbb)$ is compact, as required.

\section{Corollaries of Theorem~\ref{t_comp_KTp0}}\label{s_final}

\subsection{A family of equivalent conditions}
We introduce some notation from \cite{CCKP22TAMS}.
As usual, let $\ph^\ast$ and $\psi^\ast$ denote the boundary values on $\Tbb$ of $\ph$ and $\psi$, respectively.
Given $\ph, \psi \in\scw(\Dbb)$, $0 < p < \infty$ and $s \ge 0$, 
define
\begin{align*}
  \lambda_{p,s}(a) 
&= \lambda^{\ph, \psi}_{p,s} (a) \\
&=
\int_{(\ph^\ast)^{-1}[Q_s(a)]}
\left|
\frac{\ph^\ast - \psi^\ast}{1 - \overline{a}\psi^\ast}
\right|^p
dm +
\int_{(\psi^\ast)^{-1}[Q_s(a)]}
\left|
\frac{\ph^\ast - \psi^\ast}{1 - \overline{a}\ph^\ast}
\right|^p
dm,
\end{align*}
where
\[
Q_s(a) = \{z \in  \overline{\Dbb} : |z-a| < 2^{s+1} (1 - |a|)\}
\]
for $a\in \Dbb$.

The following corollary extends \cite[Corollary~4.1]{CCKP22TAMS}
from $1 \le p < \infty$ to all $0< p < \infty$.

\begin{cor}
Let $0 < p < \infty$ and $\ph, \psi \in \scw(\Dbb)$. The following assertions
are equivalent:
\begin{itemize}
  \item [(a)] $C_\ph-C_\psi$ is compact on $H^p(\Dbb)$;
  \item [(b)] $\lim_{|a|\to 1-} \frac{\lambda_{p,s}(a)}{1-|a|} = 0$ for some (equivalently, for all) $s\ge 0$.
\end{itemize}
\end{cor} 

\begin{proof}
As indicated above, by \cite[Corollary~4.1]{CCKP22TAMS}, (a) $\Leftrightarrow$ (b) for $1 \le p < \infty$.
Next, by \cite[Property~(4.2)]{CCKP22TAMS}, assertion~(b) does not depend on $p$ for $0 < p < \infty$.
Thus, an application of Theorem~\ref{t_comp_KTp0} completes the proof of the corollary.
\end{proof}

\subsection{Symbols $\ph$ and $\psi$ of bounded multiplicity}

\begin{cor}
Let $0 < p < \infty$ and let $\ph, \psi \in \scw(\Dbb)$ be of bounded multiplicity. Then 
$C_\ph-C_\psi$ is compact on $H^p(\Dbb)$ if and only if
\[
\lim_{|z|\to 1-} \left[\frac{1-|z|^2}{1-|\ph(z)|^2} + \frac{1-|z|^2}{1-|\psi(z)|^2}\right] 
\dist(\ph(z), \psi(z)) = 0,
\]
where $\dist$ denotes the pseudohyperbolic distance.
\end{cor} 
\begin{proof}
For $p=2$, the required equivalence holds true by \cite[Theorem~4.5]{CCKP22TAMS}.
Thus, an application of Theorem~\ref{t_comp_KTp0} finishes the proof of the corollary.
\end{proof}

\bibliographystyle{amsplain}

\end{document}